\documentclass{amsart}
\usepackage[dvips]{graphicx}
\usepackage[arrow,matrix,arc,curve]{xy}
\newtheorem{thm}{Theorem}[section]

\newtheorem{lem}[thm]{Lemma}

\theoremstyle{definition}

\theoremstyle{remark}

\numberwithin{equation}{section}

\begin{document}

\title[Framing and the Self-Linking Integral]{Framing and the Self-Linking Integral}
\author{Daniel Moskovich}
\address{Research Institute for Mathematical Sciences,
Kyoto University, Kyoto, 606-8502 JAPAN}
\email{danieru@jmail.co.jp}

 \subjclass{Knot Theory, Differential Geometry}

\date{04.03.21}
\begin{abstract}
The Gauss self-linking integral of an unframed knot is not a knot
invariant, but it can be turned into an isotopy invariant by
adding a correction term which requires adding extra structure to
the knot. We collect the different definitions/theorems/proofs
concerning this correction term, most of which are well-known (at
least as folklore) and put everything together in an accessible
format. We then show simply and elegantly how these approaches
coincide.
\end{abstract}
\maketitle

\section{Introduction}

In 1833 Carl Friedrich Gauss, in his investigation of
electromagnetic theory, discovered an integral formula for the
linking number of two space curves. If $\gamma_{0}$ and
$\gamma_{1}$ are disjoint embeddings of $S^{1}$ into $S^{3}$---
i.e two disjoint space curves, and if $\Phi:S^{1}\times
S^{1}\rightarrow S^{2}$ is the map that assigns to each $(x,y)\in
S^{1}\times S^{1}$ the unit vector from $\gamma_{0}(x)$ to
$\gamma_{1}(y)$, with $\omega$ defined as the volume form on
$S^{2}$, the Gauss integral is
$$
\frac{1}{4\pi}\int\limits_{S^{1}\times S^{1}}\Phi^{*}\omega
$$
 It is natural to ask what happens if we take $\gamma_{0}$ to equal
 $\gamma_{1}$, in other words if we want to find a knot
 invariant analogous to the linking number of a two component
 link. But in the case of a knot we run into the problem
 that $\Phi(x,x)$ is not defined (how do we define the direction from a point to itself?).
 So rather than
 $\Phi$
 being a function from $S^{1}\times S^{1}$ to
 $S^{2}$ it is instead a function from
 $C_{2}(S^{1}):=\{(x,y)\in S^{1}\times S^{1}\mid x\neq y\}$ (the configuration space of two points on
 $S^{1}$)
 to $S^{2}$. The natural way to transport the Gauss
 integral to the case of a knot is then
$$
\eta(\gamma):=\frac{1}{4\pi}\int\limits_{C_{2}(S^{1})}\Phi^{*}
\omega
$$
The problem now is that since $C_{2}(S^{1})$ is not a compact
space, we are not guaranteed that the integral converges. There
are two ways that we might try to solve this problem.\par

 \ \par
\begin{enumerate}
\item We could compactify the configuration space, and examine by
how much $\eta(\gamma)$ `fails' to be invariant. By Stokes'
Theorem, we find that this quantity depends on the boundary of the
compactified configuration space. We seek to eliminate this
boundary by pasting some extra discs $D_{0}$ and $D_{1}$ onto our
space, thus renormalizing the integral. Our invariant will then be
$\eta(\gamma)$ plus a correction term which will depend on a
``swaddling'' map $\beta$, the extension of $\Phi$ to the extra
discs  $D_{0}$ and $D_{1}$. The invariant thus constructed will
depend on the initial curve and on a choice of a \textbf{homotopy
class for the swaddling map $\beta$}. \
\par
\item We may look at the linking number of two copies of the knot,
when the second copy is ``pushed off'' to a distance of
$\varepsilon$ from the first copy, and calculate what happens as
$\varepsilon\rightarrow 0$. But this linking number will depend on
\textbf{which direction} we decided to push off the second copy of
the knot in relation to the first copy at each point, which
implicitly specifies a knot framing. As $\epsilon\rightarrow 0$,
the limit will not necessarily be an integer, forcing the
introduction of a framing dependant correction term which will
turn out to be the total holonomy of the curve. The invariant thus
constructed will depend on the initial curve, along with a choice
of \textbf{framing} for it.
\end{enumerate}

It is not at all clear at first what these two constructions
should have to do with one another. The aim of this note is to
present both approaches in a clear and accessible fashion, and to
showing how they relate in basic differential geometric terms. We
are not trying to say anything new per se, but rather to present
definitions, facts, and proofs most of which are well known, at
least as folklore, in a simple and accessible format.

\subsection{Historical remarks}
The importance of the Self-Linking Integral is that it is the most
simple and basic example of presentation of a Vassiliev invariant
as a configuration space integral. Moreover, as the work of Bott
and Taubes \cite{BT} (see also \cite{AF}) shows, this integral
plays a basic role as a correction term for anomalies in the
definition of more general finite-type ``self-linking
invariants''. In this regard, this invariant constitutes a basic
ingredient in the understanding of the Chern-Simons invariants of
knot theory.
\par

The first effective `renormalization' of the Gauss integral by
adding a correction term was carried out by Calugareanu
\cite{Cal}, and later by Pohl \cite{Poh} in the case of a closed
space curve with nowhere vanishing torsion.\par

Who the first person was to extend the invariant to curves that
may have a non-vanishing torsion I do not know. The `holonomy'
construction (the second method above) appears the more common,
and is used for instance by Polyakov \cite{Pol} (see also Tze
\cite{Tze}),  by Bott and Taubes \cite{BT} and by Bar-Natan
\cite{BN}. Meanwhile, the swaddling construction (the first method
mentioned above) is preferred by Dylan Thurston \cite{Th} and
appears more recently in papers by Poirier (\cite{Poi} and by
Lescop \cite{LeS}. The Poirier paper also gives a brief
explanation for the equivalence between the two constructions
(\cite{Poi}, remark 6.17).
\par

\subsection{Acknowledgements}
I would like to thank Prof. Dror Bar-Natan, who had the idea for
this paper, and who helped me enormously at every stage of getting
it written. I also wish to thank Tomotada Ohtsuki for his many
helpful suggestions, Emanuel Farjoun for his generous help with
the topological side of things, and Raymond Lickorish for his
helpful comments.


\section{The swaddling map construction}

Our problem when attempting to transport the Gauss self-linking
integral to knots is what to do about points of the form $(x,x)$,
for which the Gauss map is not defined. We cannot simply ignore
them, since this would force us to integrate over the space
$C_{2}(S^{1})$ which is not compact, so that the Gauss integral
would not be guaranteed to converge. If we want to find an
invariant based on the concept of self-linking, we have no choice
but to extend the Gauss function to points of the form $(x,x)$ in
some way. The problem is that as points in $C_{2}(S^{1})$ approach
points on the diagonal, the Gauss map has two limits--- the
forward and the backward sweeping tangents.

\subsection{Compactifying the configuration space}
Let us define $\overline{C_{2}(S^{1})}$, a compactification of
$C_{2}(S^{1})$, by pasting two copies of the diagonal,
$\Delta_{0}$ and $\Delta_{1}$ to its missing diagonal
$\{(x,x)|x\in S^{1}\}$, as shown in the diagram below.

\[\xy 0;/r4pc/: (0,0)="a"; (1,1)="b" **@{--}; (0,1)="c" **@{-};
?(0.45)*{>>}; (0,1)="c"; "a" **@{-}; ?(0.45)*{\wedge}; (2,0)="a";
(3,1)="b" **@{-}; (2,1)="c" **@{-}; ?(0.45)*{>>}; (2,1)="c"; "a"
**@{-}; ?(0.45)*{\wedge}; (0,-.2)="a"; (1,.8)="b" **@{--};
(1,-.2)="c" **@{-}; ?(0.45)*{\wedge}; (1,-.2)="c"; "a" **@{-};
?(0.45)*{>>}; (2,-.2)="a"; (3,.8)="b" **@{-}; (3,-.2)="c" **@{-};
?(0.45)*{\wedge}; (3,-.2)="c"; "a" **@{-}; ?(0.45)*{>>};
(.5,1.3)*=0{\textbf{$C_{2}(S^{1})$}};
(2.5,1.3)*=0{\textbf{$\overline{C_{2}(S^{1})}$}};
\endxy\]
\ \par Points on $\Delta_{0}$, which are limits of the form
$\lim_{y\rightarrow x^{+}}(x,y)$, shall be denoted $(x,x^{+})$,
with points on $\Delta_{1}$ correspondingly denoted $(x,x^{-})$.
At these boundaries of $C_{2}(S^{1})$, the Gauss map converges to
the tangent vector to the curve, sweeping either forwards or
backwards depending on whether its input converges to a point in
$\Delta_{0}$ or to a point in $\Delta_{1}$. This allows us to
solve the problem of how to extend the Gauss map to the diagonal.
$\Phi$ can be extended smoothly to a function
$\overline{\Phi_{\gamma}}: \overline{C_{2}(S^{1})}\rightarrow
S^{2}$, defined as
$\overline{\Phi_{\gamma}}(x,x^{\pm}):=\pm\dot{\gamma}$ on the
boundary.

\subsection{Checking invariance}
Let $H: S^{1}\times I\rightarrow S^{3}$ be a one parameter family
of curves. For $t \in I$, let us define
\begin{equation} \label{ps:inw}
\eta_{t}(\gamma):= \frac{1}{4\pi}\int\limits_{
\overline{C_{2}(S^{1})}\times\{t\}} \overline{\Phi^{*}}\omega
\end{equation}
Invariance of $\eta$ means that
$\eta_{0}(\gamma)=\eta_{1}(\gamma)$ for all $\gamma$. But:

\begin{equation} \label{ps:inv}
0=\int\limits_{\overline{C_{2}(S^{1})}\times I}
d\overline{\Phi_{H}^{*}}\omega=
\int\limits_{\partial(\overline{C_{2}(S^{1})}\times
I)}\overline{\Phi_{H}^{*}}\omega=
\eta_{1}-\eta_{0}+2\int\limits_{S^{1}\times I}
\overline{\Phi_{H}^{*}}\omega
\end{equation}

The first equality holds because $d$ and $\overline{\Phi_{H}^{*}}$
commute, and since $\omega$ is a 2-form defined on a 2-manifold,
$d\omega = 0$. The second equality is the Stokes theorem. The
third equality is simply the fact that $\partial
\overline{C_{2}(S^{1})} = 2S^{1}$.\par

From (2.2), we learn that $\eta_{1}=\eta_{0}$ if and only if
$\int\limits_{S^{1}\times I} \overline{\Phi_{H}^{*}}\omega =0$.
But we have no reason to assume that this would generally be the
case.

\subsection{Introducing a correction term}
Let us then ``cap off'' the cylinder $S^{1}\times I$ by pasting
two $D^{2}$'s to it, making it isomorphic to $S^{2}$, as shown in
the illustration below.

\begin{figure}[ht]
\begin{center}
\includegraphics[width=9cm,clip]{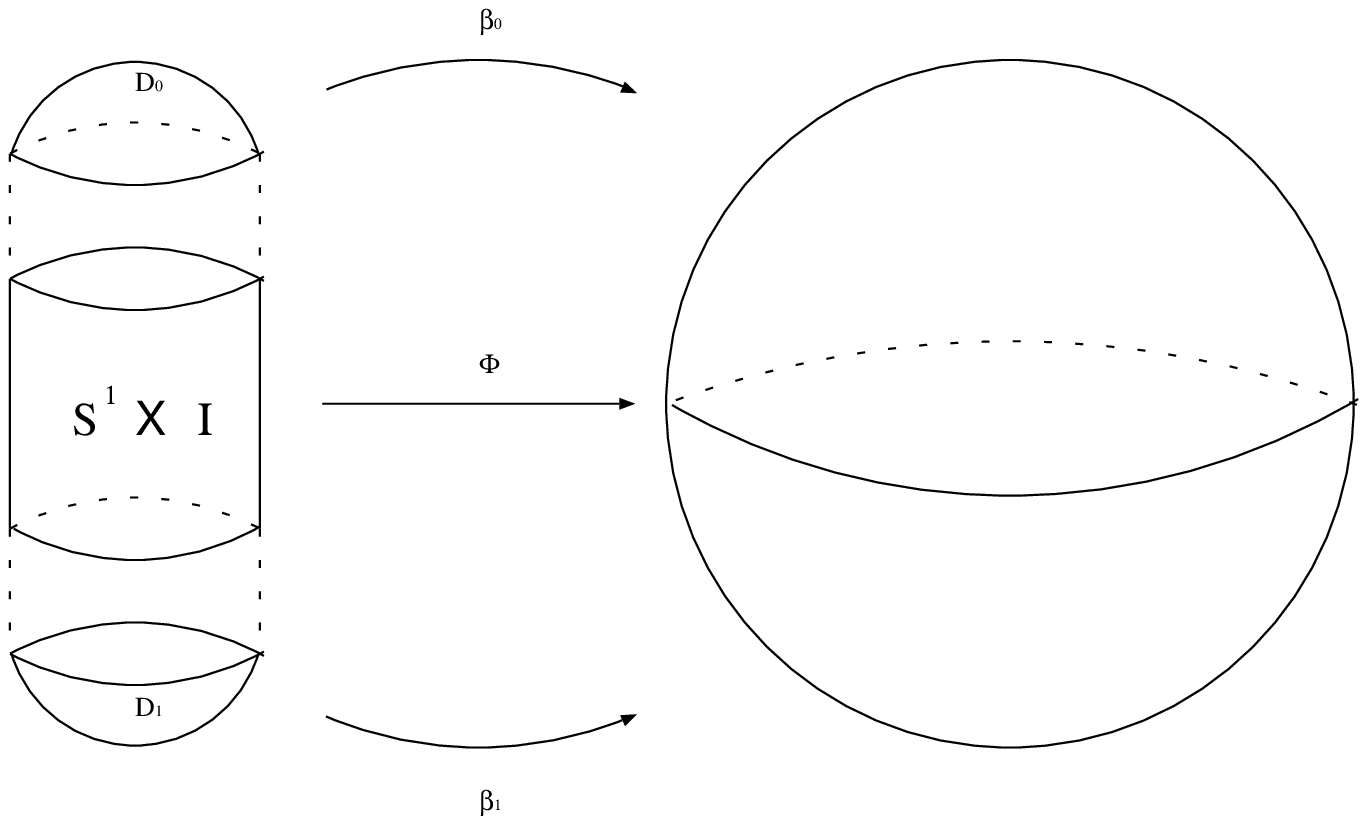}
\end{center}
\end{figure}

\ \par Let us define a ``swaddling'' map $\beta$ as a continuous
map which wraps $S^{2}$ in $D^{2}$ such that on the boundary
$\beta |_{\partial D^{2}}=\dot{\gamma}$. In our case we have two
such maps, $\beta_{0}$ and $\beta_{1}$. $\overline{\Phi_{H}}$ maps
the boundaries of $S^{1}\times I$ in ``antipodally''--- i.e. if
$x\times \{0\}\in S^{1}\times\{0\}$ maps to $y\in S^{2}$, $x\times
\{1\}\in S^{1}\times\{1\}$ maps to $-y$. Let us now define
$\beta_{0}$ as a map that maps the border of a disc to
$\pm\dot{\gamma}$, with the sign corresponding to that of
$\overline{\Phi_{H}}(S^{1}\times \{0\})$, and $\beta_{1}$ to be
$-\beta_{0}$. By abuse of notation, let us now define our total
map $\beta$ as the difference between these two maps. Our
correction term will then be defined by the equation
$$
\tau_{\beta}(\gamma):=\frac{1}{2\pi}\int_{D^{2}}\beta^{*}\omega
$$
The motivation
for this is that just like with $\gamma$, invariance of
$\tau_{\beta}$ means that $\tau_{\beta , {0}}(\gamma)=\tau_{\beta
,{1}}(\gamma)$ for all $\gamma$. But:
$$
0=\int\limits_{D^{2}\times I} d\beta^{*}\omega=
\int\limits_{\partial(D^{2}\times I)}\beta^{*}\omega= \tau_{\beta
,{1}}-\tau_{\beta ,{0}}+2\int\limits_{S^{1}\times I}
\beta^{*}\omega \label{motv}
$$
So $\tau_{\beta}$ is ``at a distance'' of
$2\int\limits_{S^{1}\times I} \beta^{*}\omega$ from being
invariant.
\par
We have found that $\partial\tau_{\beta}=\partial\eta$, proving
the following:
\begin{thm}
$\eta-\tau_{\beta}$ is an invariant of ordered pairs of a knot and
an integer specifying the homotopy class for the map $\beta$.
\end{thm}

Indeed, a simple Stokes' theorem argument shows that for two
homotopic $\beta$'s give us the same $\tau_{\beta}$ (we can push
the difference to the boundary, where the two $\beta$'s will
coincide). Moreover, our new invariant assumes integer values,
because $\beta$ wraps the disc around $S^{2}$ a whole number of
times (it has to, since $\pi_{2}(S^{2})=\mathbb{Z}$, and so
$\int_{D^{2}}\beta^{*}\omega$ assumes values in
$4\pi\mathbb{Z}$).\par Thus we find that
$sl(\gamma):=\eta(\gamma)-\tau_{\beta}(\gamma)$ (the
``self-linking number'' of $\gamma$) is an integer-valued
invariant of closed space curves along with a choice of homotopy
class of swaddling maps $\beta$. But just how much information are
we adding about the knot when we specify a homotopy class for a
swaddling map?

\subsection{Relating $\tau_{\beta}$ to the total torsion of a
space curve}

In passing, we may note that for a curve with nowhere vanishing
curvature, $\tau_{\beta}$ corresponds to the notion of the
correction term for the self-linking number as it was first
defined by Calugareanu and later by Pohl, as the total torsion of
a space curve. In section 2.4, we discovered that the correction
term $\tau_{\beta}$ is equal, modulo $\mathbb{Z}$, to the area on
$S_{2}$ covered by the map $\beta$. By the Gauss-Bonnet Theorem,
this is equal to the total curvature of $\dot{\gamma}$. But if
$\gamma$ has nowhere vanishing curvature, then this equals simply
the total torsion of the space curve $\gamma$.


\section{The holonomy construction}

In the previous section, we transported the concept of a linking
number from links to knots by compactifying the configuration
space and pasting pieces onto it in order to ``force'' the Gauss
integral to converge and to give us an invariant. There is another
way to approach the problem however. We already know that the
linking number is a link invariant--- well then, let's pretend
that our knot is a link! If we take two copies of our knot that
are only an $\varepsilon$ apart, and then see what happens when
$\varepsilon$ goes to zero, we may utilize the known invariance of
the linking number for links in order to directly conjure up an
invariant for knots.\par

In the approach we are now taking, we have one `stationary' curve,
which we shall label $\gamma_{0}$, for which we choose a smooth
framing $n(t)$. We then take the curve obtained made up of points
$t+\varepsilon\cdot n(t)$ for $t\in \gamma_{0}$, with
$0<\varepsilon\in \mathbb{R}$, which we shall denote
$\gamma_{\varepsilon}$. Now we let the `mobile' $\gamma_{1}$
descent towards the `stationary' $\gamma_{0}$. What we want to
know is what happens the self-linking integral when they
touch.\par

The classical approach here is to take the limit as $\varepsilon
\rightarrow 0$ of the Gauss integral, which involves writing and
partially calculating an explicit integral. This leads us to what
physicists call the `point-splitting regularization integral'. Tze
\cite{Tze} quotes a `simplified approach' which he credits to
Anshelevich, as quoted in an article about the twisting of strands
of DNA \cite{KaV} (!), which leads to the conclusion that the
correction term must $\frac{1}{2\pi}$ times the total holonomy of
the curve. But the same result may be obtained in a more
elementary way by making use of a technique we have already
used--- that of a `swaddling map'. This will also help us to
visualize why and when the two constructions for the error term of
the self-linking integral will coincide.

 \subsection{Two Ways of Looking at the Same Thing}

\subsubsection {Don't take limits- compactify!}

Rather than thinking of limits of integrals, let's compactify the
space of pairs $$ \mathbb{L}^{\overrightarrow{n}}(S^{1}):=\{(t, t+
\varepsilon\cdot n(t))\mid t\in \text{\{a closed space curve\}},
0<\varepsilon\in\mathbb{R}\}$$ by pasting something onto the
boundary $\mathbb{L}^{\overrightarrow{n}}_{0}(S^{1}):=
\overline{\mathbb{L}}^{\overrightarrow{n}}(S^{1})\mid_{\varepsilon=0}$,
where the overline denotes topological closure. Our problem, as
usual, is that $\Phi(x,x)$ is not defined, and what this second
approach gives is a way of defining it via a limit which keeps
track of the information which is relevant to the Gauss
integral--- direction--- and thus tells us what it wants
$\Phi(x,x)$ to be.
\par

Thus, the space we must paste on should consist of pairs
$(x,\theta)$ in which we store the ``address'' $x$ of the point,
as well as the direction from which $x$ is coming in on
$\gamma_{\varepsilon}$. We may depict the newly created boundary
of our compactified space as a continuum of pieces which can be
schematically depicted like this:

\[\xy
0;/r3pc/: (2,0)="a"; 0 **@{-}; (0,.1)="b"; (.8,.1)="c" **@{-};
(1.2,.1)="d" **\crv{(1,.4)}; (2,.1)="f" **@{-}
\endxy\]

The leftmost point of the semicircle corresponds to a point
$\gamma_{\varepsilon}(t)$ coming in to $\gamma_{t}$ on a backward
sweeping tangent, the rightmost one corresponds to the point
coming in on a forward sweeping tangent, and the apex corresponds
to the point coming in straight off the normal. The tangent and
the normal at each point define the framing for the knot, so that
we see that we have not lost any information. Let us define
\begin{equation}
\tau_{\phi}:=\frac{1}{2\pi}\int_{\mathbb{L}^{\overrightarrow{n}}(S^{1})_{0}}\overline{\Phi^{*}_{\gamma}}\omega
\end{equation}

Our main claim in this section is:
\begin{thm}
$\eta+\tau_{\phi}$ is an invariant of framed knots.
\end{thm}

To see that $\eta+\tau_{\phi}$ is in fact the limit of the
integral of the pullback of the volume form, note that the domain
at every $\varepsilon >0$ where $\varepsilon$ is fixed is
cobordant to the integral at `the bottom', where $\varepsilon=0$,
and hence via a standard Stokes' theorem argument, the
``point-splitting regularization integral'' used by \cite{Tze} is
the same as the integral along
$\mathbb{L}^{\overrightarrow{n}}_{0}(S^{1})$.\par

\subsubsection {``$C$-swaddling''}

Let us eliminate the boundary of our configuration space this time
in a different way--- rather than pasting two discs onto the
$S^{1}\times I$ boundary, let us instead paste another
$S^{1}\times I$ to it (which we shall call $C$ for cylinder),
completing the cylinder $C_{2}(S^{1})$ to a torus
$C_{2}(S^{1})\cup_{\partial}C$ as shown in the illustration below:
\
\par
\[\xy
0;/r3pc/: (-2,1)="a"; (-1,1)="b" **\crv{(-1.5,1.2)}
**\crv{(-1.5,0.8)}, (-2,0)="c" **@{-}; (-1,0)="d"
**\crv{~*=<3pt>{.}(-1.5,.2)} **\crv{(-1.5,-.2)}; "b" **@{-};
(0,1)="e"; (1,1)="f" **\crv{(.5,1.2)} **\crv{(.5,0.8)}, (0,0)="g"
**@{-}; (1,0)="h" **\crv{~*=<3pt>{.}(.5,.2)} **\crv{(.5,-.2)}; "f"
**@{-}; "a" **\crv{~*=<5pt>{.}(-.5, 2.1)};
 (-1,1)="b"; "e"
**\crv{~*=<5pt>{.}(-.5,1.3)}; (-2,0)="c"; "h"
**\crv{~*=<5pt>{.}(-.5,-1.1)}; (-1,0)="d"; "g"
**\crv{~*=<5pt>{.}(-.5,-.3)} {\ar@{}|{C_{2}(S^{1})}"a";(-1,0)}
{\ar@{}|{C}"e";(1,0)}
\endxy\]
\ \par

Then we define a new swaddling map $\phi:C\rightarrow S^{2}$ which
maps the boundaries of $C$ to the curves $\pm\dot{\gamma}$ in such
a way that $\overline{\Phi_{H}}$ and $\phi$ combine to give us a
continuous map $\overline{C_{2}(S^{1})}\cup_{\partial}C\rightarrow
S^{2}$.
\par
As the title to this section suggests, we would like to show that
in 3.1.1 and in 3.1.2 we have done one and the same thing (modulo
$4\pi$)--- that in point of fact, there is no difference between
compactifying the space as in section 3.1.1. by storing framing
information on the boundary, and between eliminating the boundary
of $\overline{C_{2}(S^{1})}$ by the ``$C$-swaddling'' method as we
have done in this section.\par

\subsubsection{Every framing gives a $C$-swaddling}

The integral on the bottom splits into two parts--- the ``normal''
pieces in which we have just the standard Gauss self-linking
integral, and the ``bumps''. On the bumps, $\theta$ goes from the
tangent to the normal to minus the tangent--- i.e. it traces out
the image of a line between two points of the form $\{x\}\times
\{1\}$, $\{x\}\times \{0\}$ on $S^{1}\times I$, thus defining a
$C$-swaddling, as the path can be assumed to be a great circle for
all point on the knot.

We thus see that the ``$C$-swaddling'' construction gives us
$\tau_{\phi}$, in just the same way as the construction in the
first section gave us $\tau_{\beta}$.

\subsubsection{Every $C$-swaddling gives a framing}

The fact that every $C$-swaddling gives a framing follows from the
following lemma.

\begin{lem} \label{E:matz}
Every $C$-swaddling map is homotopic via $C$-swaddling maps to a
map in which the path on $S^{2}$ from $\dot{\gamma}(s)$ to
$-\dot{\gamma}(s)$ for each $s$ in the knot is a great circle
between the two points.
\end{lem}

\begin{proof}

\begin{enumerate}
    \item[Step 1]
    Every $C$-swaddling map is homotopic to a boundary-fixing
    diffeomorphisms from the cylinder to itself composed with a
    mapping from the cylinder to $S^{2}$ taking the boundaries of
    the cylinder to $\dot{\gamma}$ and to $-\dot{\gamma}$ and
    mapping the line $s\times\{i\}$, $i\in \mathbf{I}$ to a
    great circle from  $\dot{\gamma}(s)$ to
    $-\dot{\gamma}(s)$ for each $s\in S^{1}$. Let us fix the
    second map in this composition and call it $g$.
    By \cite{Gra} there is a $\mathbb{Z}$-worth of
    boundary-preserving diffeomorphisms of the cylinder
    $S^{1}\times \mathbf{I}$ to itself up to homotopy by such
    diffeomorphisms. The generator of the this homotopy group is a Dehn twist
    about a boundary-parallel curve.

    \item[Step 2] Let $p_{s}$ be the image of $s\times {I}$ in the
    $S^{2}$ under a power of the generator which we found in step
    1. Let $D_{\varepsilon}$ be a disc in $S^{2}$ of radius
    $\varepsilon$ centred at the image of $s\times \{1\}$. The
    homotopy which fixes $s\times \{0\}$ and $s\times\{1\}$ and
    revolves $p_{s}\cap\partial D_{\varepsilon}$ in a full circle
    around $\partial D_{\varepsilon}$ for all $s\in S^{1}$ undoes
    the Dehn twist.\par

    Concretely, let's look at $\varphi_{t}$ $t\in \mathbf{I}$
    which takes $p_{s}$ to a smoothing of the curve obtained by
    taking a radius of $D_{\varepsilon}$ at angle $2\pi t$, then travelling around
    the circumference of $D_{\varepsilon}$ clockwise until we hit $p_{s}$
    and continuing with $p_{s}$. $\varphi_{1}(p_{s})$ can be taken by homotopy to
    a curve in the image of the cylinder on $S^{2}$, which is
    $p_{s}$ with a Dehn twist added.

    \item[Step 3] Thus for each mapping $g$ considered in step 1,
     the map from the cylinder has a single representative $g$ composed with the identity.
     There are a $\mathbb{Z}$-worth of choices of $g$, corresponding to taking the path on
     the cylinder from $s\times {i}$, $i\in \mathbf{I}$, to the great circle which wraps around the
     circle $n$ times for all $n\in \mathbb{Z}$.
    \end{enumerate}
\end{proof}

Choosing the midpoint for each such great circle gives us a normal
to the knot at $s$. In other words, given a family of
forward-sweeping tangents to the knot, a $C$-swaddling map gives
us a smooth family of normals to the knot, thus giving the knot a
framing.

\subsection{Holonomy through $C$-swaddling}
To show that the correction term $\tau_{\phi}$ that we get is the
total holonomy, we must first represent $\phi^{*}\omega$ as the
pullback of an element of $SO(3)$. For this purpose, as the $\phi$
swaddling map is a smooth extension of the Gauss map, let us
redefine $\overline{\Phi_{\gamma}}(x,y)$ to be
$\Phi_{\gamma}(x,y)$ when $x\neq y$ and $\phi$ on the boundary.

\subsubsection{Transporting the pullback of the volume form to
SO(3)} Since we are now moving into SO(3), we shall convert the
discussion into the language of framings. Let us break
$\overline{\Phi_{\gamma}}(x,y)$ into the mapping $\phi$ from
$S^{1}$ to $SO(3)$ composed on a mapping $e_{1}(x,y)$ from $SO(3)$
to $S^{2}$ . Following \cite{Poh}, let $e_{2}$ be the unit vector
normal to $e_{1}$ on the plane spanned by $e_{1}$ and the tangent,
extending smoothly to the normal defined by $\phi$ on the
boundary. We shall then define $e_{3}$ to be $e_{1}\times e_{2}$
at every point. The following lemma is due to William Pohl
\cite{Poh}.

\begin{lem}
$$
e_{1}^{*}\omega=d(de_{3}\wedge e_{2})
$$
\end{lem}

\begin{proof}
For every $1\leq i_{3}$ in $\mathbb{N}$, $e_{i}$ defines a
function $x_{i}$ by means of the relation $e_{i}^{0}\cdot
v=x_{i}(v)$ for any vector $v\in\mathbb{R}^{3}$, when $e_{i}^{0}$
denoted the restriction of $e_{i}$ to the point
$\overrightarrow{0}$. The volume form in $\mathbb{R}^{3}$ is then
given by the expression $$x_{1}(dx_{2}\wedge dx_{3})+ \text {
cyclic permutations.}$$

Let us pull back the volume form via $e_{1}$. $e_{i}^{*}x_{j}=
e_{i}\cdot e_{j}^{0} = \delta_{i,j}$, therefore
$$e_{1}^{*}\omega= d(e_{1}^{*}x_{2})\wedge d(e_{1}^{*}x_{3})=de_{1}\cdot e_{2}^{0}
\wedge de_{1}\cdot e_{3}^{0}$$ by Leibnitz's rule. There was
nothing special about our choice of 0 as the point by which to
define the functions $x_{i}$, therefore we have $de_{1}\cdot
e_{2}\wedge de_{1}\cdot e_{3}$. \par

But $de_{i}\cdot e_{i}=0$, and by differentiating this equality we
find that $de_{i}\cdot e_{j}=-de_{j}\cdot e_{i}$ and so this
equals
\begin{equation}
\label{de:wej} (de_{3}\cdot
e_{1})\wedge(de_{2}\cdot e_{1})
\end{equation}

Now we remember that the $e^{i}$'s are an orthonormal to one
another, and therefore they satisfy the equality
$$
de_{3}\cdot de_{2}=\sum_{i=1}^{3}(de_{3}\cdot
e_{i})\wedge(de_{2}\cdot e_{i}) \label{wedg}
$$

But according to (3.1) this is exactly $e_{1}^{*}\omega$, and so
we have found that
$$
e_{1}^{*}\omega=de_{3}\cdot de_{2}=de_{3}\cdot de_{2}+
d^{2}e_{3}\cdot e_{2}=d(de_{3}\wedge e_{2}) \label{wege}
$$
\end{proof}

\subsubsection{Relating $\tau_{\phi}$ to the total torsion of a
space curve} In $\gamma$ has nowhere vanishing curvature, we can
use Lemma (3.2) to show that this correction term as well is equal
to the total torsion of the curve $\gamma$. Here $e_{1}$ is the
tangent, $e_{2}$ the normal, and $e_{3}$ the binormal, so
$$de_{3}\cdot e_{2}=db\cdot n$$ $-\tau\cdot n\cdot n=-\tau$, so by
the Frenet equations, $b'\cdot n=-\tau$ so $db\cdot n=-\tau\cdot
ds$.

Thus, we see that for a curve with a nowhere-vanishing curvature,
$$\tau_{\phi}=\int_{S^{1}}\tau ds \label{crv}$$
which is again the total torsion of the space curve $\gamma$.

\subsubsection{Making sense of it all}
The last step of our argument is just the Stokes' theorem. The
domain of $\overline{\Phi_{\gamma}}$ is our ``cylinder
compactification'' of $C_{2}(S^{1})$. Pulling back the volume form
via this map, when restricted to $C$, will then by Stokes' Theorem
be equivalent to pulling back $de_{3}\wedge e_{2}$ via $\phi$
along $C$'s boundary. But here $\phi$ gives us the tangent,
$e_{2}$ the normal, and $e_{3}$ the binormal,
$\phi^{*}de_{3}\wedge e_{2}$ is the triple product
$(\dot{\gamma},n,\dot{n})$. Here though we have
$\int_{S^{1}}(\dot{\gamma},n,\dot{n})ds=
\int_{S^{1}}(\dot{\gamma},n,(n(s)+\dot{n}(s)ds))
=\int_{S^{1}}(\dot{\gamma},n,\dot{n}(s+ds))$. The last integral is
measuring ``by how far'' the normal has strayed from its initial
position at $t=0$ by the time we get to $t=L$. In other words,
$\tau_{\phi}$ is measuring the \textbf{total holonomy} of the
curve $\gamma$, with respect to the Reimannian connection on the
normal bundle to the curve.

\subsection{Equivalence to total torsion (again)}

 Here again we
have an easy proof that the correction term of the self-linking
integral equals the total torsion. When our curve has a Frenet
frame $(t,n,b)$ with curvature $\kappa$ and torsion $\tau$, the
Frenet equations give us
$$
\tau_{\phi}=\frac{1}{2\pi}\int_{S^{1}}(\dot{\gamma},n,\dot{n})ds=
\frac{1}{2\pi}\int_{S^{1}}(\dot{\gamma},n,(\tau b-\kappa
t))ds=\frac{1}{2\pi}\int_{S^{1}}\tau ds \label{tort}
$$


\section{Equivalence of the two constructions}

In the previous sections, we have presented two alternative ways
of introducing a correction term to the Gauss self-linking
integral for a knot, making it an invariant. We do not know yet
whether these two methods are equivalent, and there is no reason a
priori to assume that this should be the case. Why should choosing
a homotopy class for a swaddling map have anything to do with
choosing a framing for a knot? In both of these approaches, we
reach the image on $S^{2}$ via the Gauss map of the tangent bundle
to an embedding into $S^{3}$ of $S^{1}$, but in the first approach
we come to this image by first embedding $S^{1}$ into $D^{2}$ and
then getting to $S^{2}$ via the swaddling map $\beta$, while in
the second approach we first map to $SO(3)$ by choosing a framing
(we shall call this map $\phi$), and then map down from there onto
$S^{2}$. The situation is schematically depicted in the
commutative diagram below: \
\par
\[
\xymatrix{D^2 \ar[drr]_(.6)\beta
&& SO(3) \ar[d]^{e_{1}} \\
S^{1} \ar@{^{(}->}[u] \ar[rr]_{\dot{\gamma}} \ar[urr]^(.6)\phi &&
S^{2}}\]

In this notation,
$\tau_{\beta}(\gamma)=\frac{1}{2\pi}\int_{D^{2}}\beta^{*}\omega$,
while $\tau_{\phi}(\gamma)=\frac{1}{2\pi}\int_{S^{1}}\phi^{*}\tau$
for $\tau$ a pre-image of omega via the map $e_{1}$. Equality of
these terms would follow from the existence of a map $\sigma$ such
that the following diagram commutes:

\[
\xymatrix{D^2 \ar[drr]_(.6)\beta \ar@{-->}[rr]^\sigma
&& SO(3) \ar[d]^{e_{1}} \\
S^{1} \ar@{^{(}->}[u] \ar[rr]_{\dot{\gamma}}  \ar[urr]^(.6)\phi &&
S^{2}}\] For in that case
$$
\tau_{\beta}(\gamma)=\frac{1}{2\pi}\int_{D^{2}}\beta^{*}\omega=\frac{1}{2\pi}\int_{D^{2}}\sigma^{*}d\tau=
\frac{1}{2\pi}\int_{S^{1}}\phi^{*}\tau=\tau_{\phi}(\gamma)
$$
The second equality stems from the fact that the diagram is
commutative, and the fourth we have already shown. So all the
action takes place around the middle equality. We have shown that
when $\beta^{*}\omega$ is transported to SO(3), it becomes the $d$
of something. So we may use Stokes' theorem to go from left to
right.\par

We can also see this easily from the swaddling map construction-
let us choose a $\beta$ mapping, pasting two discs onto the
boundaries of $C_{2}(S^{1})$, making it a compact space. Let us
choose our discs such that corresponding points on $D_{1}$ and on
$D_{0}$ map to antipodal points on $S^{2}$ via $\beta$. Cutting
out a small neighbourhood of the centres of the discs, we may glue
a cylinder between them, connecting them into a shape isomorphic
to the cylinder on which our $\phi$ map was defined. Now every
$\beta$ map can be smoothly extended to a $\sigma$ map, because
the two discs with the narrow tube connecting them is
homotopically a cylinder. \par

But as Tahl Novik observed, going from right to left in this set
of equalities we have to watch out, because
$\pi_{1}(SO(3))=\mathbb{Z}/2$, and for a path belonging to the
non-trivial homotopy class of $SO(3)$, there can exist no
pre-image via a $\sigma$ mapping.
\par

Let us note that the cylinder of the $C$-swaddling construction
can be `cut' into 2 discs if and only if it is homotopic to a
cylinder of which the `middle circle' is constant- in other words
as a framing it is homotopic to the constant framing. Then and
only then can we `pinch closed' that sphere, turning the cylinder
into two discs tangent at a point without loss of information. For
elements of the non-trivial homotopy class, this is by definition
going to be impossible. Notice that by `pinching' the cylinder
into discs, we are separating the backward sweeping tangents and
the forward sweeping tangents, which is impossible in the
non-trivial homotopy class in which these two families of tangent
vectors are one and the same.\par

But for elements of the trivial homotopy class of $SO(3)$, no such
difficulty arises. Stokes' theorem takes us from
$\tau_{\phi}(\gamma)$ to $\tau_{\beta}(\gamma)$. We have proved
then the following theorem:

\begin{thm}
$\tau_{\beta}+4\pi=\tau_{\phi}$
\end{thm}

The ``$+4\pi$'' correction is an idea of Tomotada Ohtsuki's, to
remind us that for `minimal' representatives of $\phi$ and $\beta$
the area $\phi$ covers on $S^{2}$ with the cylinder $C$ (in this
case `minimal' would be taken to mean that each ``vertical'' line
between the boundaries of the cylinder is mapped to the minimal
length line between the tangent and minus the tangent on $S^{2}$,
with appropriate sign) is the entire ball, plus the area $\beta$
covers with the two discs (and here `minimal' means simply the
minimal such positive area). In any event, modulo $4\pi$ the
correction terms are equal. The isomorphism between the two
correction terms means that in a very real sense choosing a
swaddling map $\beta$ along with the homotopy class in which it
sits is exactly the same thing as choosing a framing that is
null-homotopic as an element of $SO(3)$. We have shown then that
$\beta$ can be lifted to $\sigma$, but only for `half' our
possible choices of $\phi$.

\subsection{So which ``half'' is it?}

We have shown then that for framings which give us an element of
the trivial homotopy class of $SO(3)$,
$\tau_{\phi}(\gamma)=\tau_{\beta}(\gamma)$. We have yet to show
what framings those are.\par

Let us recall the mappings defined in section 2.3.
$\overline{\Phi_{H}}$ mapped $S^{1}\times I$ to $S^{2}$, sending
the two components of the boundary to the tangent bundle of the
knot $\gamma$ in antipodal ways, while $\beta_{0}$ and $\beta_{1}$
took the boundary of a disc, and mapped it to $\dot{\gamma}$ and
to $-\dot{\gamma}$ correspondingly. Thus, we have a
$\mathbb{Z}_{2}$ action on $\overline{C_{2}(S^{1})}$, whose action
is to flip: $(x,y)\rightarrow (y,x)$. $\overline{\Phi_{H}}$ then
descends to the quotient
$$ \label{sim}
\Phi_{H}:\overline{C_{2}(S^{1})}/\mathbb{Z}_{2}\rightarrow
S^{2}/\mathbb{Z}_{2}\simeq \mathbb{R}P^{2}
$$
But $\overline{C_{2}(S^{1})}/\mathbb{Z}_{2}$ is also just
$\mathbb{R}P^{2}$, so $\Phi_{H}$ is in fact a map from
$\mathbb{R}P^{2}$ to itself. \par $\mathbb{R}P^{2}$ is a
non-orientable space, therefore only the degree of
$\overline{\Phi_{\gamma}}$ is only defined mod 2. But the flipping
action $(x,y)\rightarrow (y,x)$ is precisely the non-trivial path
in $\mathbb{R}P^{2}$, hence the degree of
$\overline{\Phi_{\gamma}}$ must be 1 (this follows from the
topological assertion that the degree of a map is $\pi_{1}$ of
that map).\par

This gives us a complete characterization of the framings for
which $\phi$ lifts to $\sigma$- they are exactly those framings
for which the mod 2 degree of the `extended Gauss mapping'
$\overline{\Phi_{\gamma}}$ is 1. This leads us to the rather
startling conclusion that, given the blackboard framing, our
`trivial' knot turns out not to be the circle at all, but rather
the boundary of the Moebius band.\par

\section{A combinatorial description of the invariant}
By adding a correction term, we have shown that the Gauss
self-linking integral can be made to be an invariant of framed
knots. It so happens \cite{BF} that this invariant coincides with
the so-called `writhing number' of the curve, obtained by taking
the number of positive crossings and subtracting the number of
negative crossings. Thus, we have constructed an invariant
analogous to the \emph{linking number} of two disjoint space
curves.\par

There is also another combinatorial description of our invariant,
which is to my mind more appealing \cite{Poh}, \cite{Aic}. Let us
imagine the knot as a roller-coaster, with us sitting in a car
facing forwards. At every point, the tracks face away from the
knot in the direction of the normal vector, and as the car travels
along the rails, our head is always pointing ``up''. Let us also
assume that our head is locked in place, such that we can only
look straight ahead (in the direction of the tangent).\par

The roller coaster starts up, and we start moving along the track.
The car rises and falls, twists and loops, swooshing along. Every
now and again, we may see another portion of track coming up
directly into our field of vision--- Pohl calls such points
`cross-tangents'. We count these with appropriate sign, depending
on the orientation of the tracks (which way the car has gone down
them or will go down them, and the direction in which we are
currently travelling). The roller coaster stops when we return to
our initial point, and we sum up all the cross-tangents, with
appropriate signs. And we get what we have calculated in this
paper--- the Gauss self-linking integral with the appropriate
correction term, determined by the direction the tracks faced away
from the knot at each point.


\bibliographystyle{amsplain}

\end{document}